\documentclass[12pt, a4paper]{article}

\usepackage{amsmath}
\usepackage{amsthm}
\usepackage{amssymb}
\usepackage{color}
\usepackage{pdfpages}
\usepackage{url}
\usepackage{comment}
\usepackage{mathrsfs}

\usepackage[unicode,colorlinks=false,linkcolor=red,citecolor=black,urlcolor=green,pagebackref=false]{hyperref}

\theoremstyle{definition}
\newtheorem{theorem}{Theorem}
\theoremstyle{definition}
\newtheorem{proposition}[theorem]{Proposition}
\theoremstyle{definition}
\newtheorem{lemma}[theorem]{Lemma}
\theoremstyle{definition}
\newtheorem{definition}[theorem]{Definition}
\theoremstyle{definition}
\newtheorem{example}[theorem]{Example}
\theoremstyle{definition}
\newtheorem{corollary}[theorem]{Corollary}
\theoremstyle{definition}

\theoremstyle{definition}
\newtheorem{assumption}[theorem]{Assumption}
\theoremstyle{definition}
\newtheorem{remark}[theorem]{Remark}

\newcommand{\N}{\mathbb{N}}	
\newcommand{\R}{\mathbb{R}} 
\newcommand{\E}[1]{\mathbb{E}\left[#1\right]}
\newcommand{\Ebig}[1]{\mathbb{E}\bigg[#1\bigg]}
\renewcommand{\P}{\mathbb{P}}

\newcommand{\F}{\mathscr{F}}
\newcommand{\X}{\mathbb{\widehat X}}
\newcommand{\ran}[2]{\langle #1, #2 \rangle}

\newcommand{\bfsig}{\boldsymbol{\sigma}}

\newcommand{\1}{\mathbf{1}}

\newcommand{\dd}{\mathrm{d}}

\title{On the integrability of the supremum of stochastic volatility models
and other martingales}

\author{
Stefan Gerhold \\
TU Wien \\
\tt{sgerhold@fam.tuwien.ac.at}
\and
Julian Pachschw\"oll \\
Universit\"at Wien \\
\tt{julian.pachschwoell@univie.ac.at}
\and
Johannes Ruf\\
CUHK Shenzhen \& LSE\\
\tt{j.ruf@lse.ac.uk}}

\date{\today}

\numberwithin{equation}{section}
\numberwithin{theorem}{section}

\begin{document}

\maketitle

\begin{abstract}
We propose a method to bound the expectation of the supremum of the
price process in stochastic volatility models. It can be applied, for example, to the
rough Bergomi model, avoiding the need to discuss finiteness of higher moments.
Our motivation stems from the theory of American option pricing, as an integrable
supremum implies the existence of an optimal stopping time for any linearly bounded
payoff. Moreover,
we survey the literature on martingales with non-integrable supremum,
and give a new construction that yields uniformly integrable martingales with
this property. 
\end{abstract}

\smallskip

\section{Introduction}

This paper deals mainly with the question whether a nonnegative martingale $S=(S_t)_{0\leq t\leq T}$,
where $T\in(0,\infty]$ is deterministic, satisfies
\begin{equation}\label{eq:sup}
  \Ebig{\sup_{t\in[0,T]} S_t} <\infty.
\end{equation}
{}From the viewpoint of mathematical finance, this property is of interest in
the theory of American option pricing. If $S$ is the price process of the underlying
and the expectation is under the chosen pricing measure, condition~\eqref{eq:sup}
implies the existence of an optimal exercise time for all American payoffs
of at most linear growth~\cite[Theorem~D.12]{KaSh98}.
A standard example is the straddle~\cite{AlMa02},
with payoff $|S_t-K|$. In the European case, its price decomposes into a call price plus a put price, where both put and call have  strike~$K$. Under American exercise rights, such a decomposition does not hold. 
Thus,~\eqref{eq:sup}
serves as theoretical basis of any convergence analysis regarding
approximate optimal exercise strategies, as it guarantees the existence of
the limiting object.

If $T<\infty$ and $\mathbb{E}[S_T^p]<\infty$ for some $p>1$, which may depend on~$T$, then~\eqref{eq:sup} follows from Doob's
$L^p$ inequality. A convenient way of verifying finiteness of moments $\mathbb{E}[S_T^p]$ is to study the domain of the characteristic
function, and so the existence of higher moments is well understood for affine models~\cite{Ke11} and affine Volterra models~\cite{GeGePi19,KeMa20}. Further examples of stochastic volatility models with known characteristic function
can be found in~\cite[Chapter~4]{Gu12}.
The $3/2$ model admits a particularly simple statement in this regard, which has
not been made explicit in the literature: If the asset and variance processes are negatively correlated, any moment $\mathbb{E}[S_T^p]$ with exponent $p\in[1,u_+)$
is finite, irrespective of the value of~$T$. We refer to~\cite{GePi20} for the definition
of~$u_+>1$ in terms of the model parameters, and for refined tail asymptotics of the marginal
density.
For models that do not feature an explicit characteristic function, a tool that could be used to establish~\eqref{eq:sup}
is~\cite[Theorem~4.1]{Pe96}, which gives sufficient conditions for the finiteness
of the exponential Orlicz norm of the maximum of a continuous local martingale.

For rough volatility models outside the affine Volterra class, it seems that there are only results
on the \emph{non}-existence of moments so far~\cite{Ga19,Gu20}, which
motivates studying~\eqref{eq:sup} directly.
We introduce a class of generic rough stochastic volatility models for which we give a sufficient condition for~\eqref{eq:sup} based on Doob's $L^1$ inequality and a change of measure. 
We demonstrate the applicability of our approach by treating several specific instances of our model class.

We give a full characterization of~\eqref{eq:sup} for the 
multi-factor rough Bergomi model~\cite{lacombe2020asymptotics}, which contains the classical rough Bergomi model as a special case.
We further prove that~\eqref{eq:sup} holds generically for models satisfying certain linear growth conditions, effectively recovering the integrability of the supremum for volatilities given by affine Volterra equations without requiring existence of higher order moments. 
Finally, we slightly extend results in \cite{jaber2025martingalepropertymomentexplosions} by showing integrability of the supremum for multi-factor signature volatility models if and only if the price process is a true martingale.

We note in passing that the stronger condition that the supremum of the expectation in~\eqref{eq:sup} over all equivalent martingale measures be finite
is of interest as well. Under this assumption, \emph{all} arbitrage free price
processes of a linearly bounded American payoff can be represented by
Snell envelopes with respect to some equivalent martingale measure. This follows  from \cite[Theorem~13.2.9]{KaKu05}, by specializing from game
options to American options. The statement is not found in many textbooks, while
the European counterpart of the theorem, viz.\ that all arbitrage free price
processes arise from equivalent martingale measures,
is standard material in mathematical finance courses -- see~\cite[Theorem~5.29]{FoSc16}
for this (the European) statement in discrete time, which avoids technicalities.

A-priori, it is not obvious whether there actually are martingales that do \emph{not} satisfy~\eqref{eq:sup}. 
We discuss this question in Section~\ref{se:sup}.
We survey the literature and show how known results easily yield explicit examples of such martingales.
Moreover, we present a new construction that gives plenty of examples, by appropriately stopping a strict local martingale.

In terms of notation, for a given probability space, the class  $L\log L$  consists of the random variables~$X$ with 
$\mathbb{E}[|X\log |X||]<\infty$.
Throughout this paper, $d,m \in \N$ denote two integers. We also consider the lower triangle
\[
    \Delta = \{(t,s) \in [0,\infty)^2: t > s\}.
\]

\section{Integrability of the supremum in stochastic volatility models} \label{se:int}

The following lemma states a version of Doob's $L^1$-inequality~\cite{ac12,Gi86}.

\begin{lemma} \label{lem:exp_sup_est}
    Let $\left(X_t\right)_{t \geq 0}$ be a nonnegative right-continuous submartingale and let $T > 0$. Then we have 
    \begin{align*}
        \Ebig{\sup_{t \in [0,T]} X_t} 
		\leq \frac{e}{e - 1} \Big(\E{X_T \log(X_T)} + \E{X_0 (1 - \log(X_0))}\Big).
    \end{align*}
    If the submartingale is additionally uniformly integrable, the extended process $\left(X_t\right)_{t \in [0,\infty]}$ is a submartingale as well by Doob's $L^1$-convergence theorem. In this case, the inequality also holds for $T=\infty$. 
\end{lemma}

Doob's original statement has $(\log(\cdot))^+$ instead
of $\log$, which is slightly less convenient for
our application.
While some sources present the inequality in discrete
time, it is clear that we can easily pass to continuous time, by means of the monotone convergence theorem.
We note in passing that an application of the inequality to pathwise hedging is presented in~\cite[Section~5]{AcBePeSc16}.

We start with a rather generic modeling setup for the spot price process~$S$ that includes various rough and non-rough stochastic volatility models as special cases.
Our main theoretical results on the integrability of the supremum are given below in Theorem~\ref{thm:int_of_sup} and Corollary~\ref{thm:2:int_of_sup}. We show in the next section that various commonly used models are in scope of our framework.

\begin{definition}[Progressive functional]
    For $t \geq 0$ define the coordinate projection $\pi_t : C([0,\infty),\R^m) \to \R^m$ by $\pi_t(w) = w(t)$. Define the sigma algebra $\mathcal{G}_t = \sigma(\pi_s: s \in [0,t])$ and the filtration $\mathbb{G} = (\mathcal{G}_t)_{t \geq 0}$. Then a function $g: [0,\infty) \times C([0,\infty), \R^m) \to \R^d$ is a progressive functional if for each $t \geq 0$, the restriction
    \begin{align*}
        g: [0,t] \times C([0,\infty),\R^m) \to \R^d
    \end{align*}
    is $\mathcal{B}|_{[0,t]} \otimes \mathcal{G}_t$-measurable.
\end{definition}
 A progressive functional depends at each time $t \geq 0$ only on path information up until time $t$. In particular, for an adapted and continuous $\R^m$-valued process $X$ and a progressive functional $g: [0,\infty) \times C([0,\infty),\R^m) \to \R^d$, the map $(\omega,t) \to f(t,X(\omega))$ is progressive in the usual sense.

\begin{definition}[Generic stochastic volatility model] \label{def:generic_sv_model}
    Set $E = [0,\infty) \times C([0,\infty),\R^m)$ and consider progressive functionals
    \begin{align*}
        &\psi: E \to \R, \quad \mu: E \to \R^m, \quad \sigma^S: E \to \R^{1 \times d} , \quad \sigma^Y: E \to \R^{m \times d} 
    \end{align*}
    with $\|\sigma^S\|_2  = 1$,
    and some measurable kernel $K: \Delta \to [0,\infty)$. 
    We consider a filtered probability space $(\Omega, \F, (\F_t)_{t \geq 0}, \mathbb{P})$ , where $\F_0$ is $\P$-trivial, satisfying the usual conditions and supporting an $\R^{d}$-valued Brownian motion $W$.
    Then a \textit{generic stochastic volatility model} is a continuous and adapted $\R^{1+m}$-valued process $(S,Y)$ on that space such that, for any $t \geq 0$,
    \begin{multline*}
        \int_0^t \psi(s,Y)^2 \dd s 
        +
        \int_0^t K(t,s) |\mu(s,Y)| \dd s \\
        +
        \int_0^t K^2(t,s)
        \mathrm{tr}(\sigma^Y(s,Y)\sigma^Y(s,Y)^\top) \dd s < \infty
    \end{multline*}
    and 
    \begin{align}
        S_t
        &= 
        S_0
        \exp
        \bigg(
            \int_0^t \psi(s,Y) \sigma^S(s, Y) \dd W_s - \frac12 \int_0^t \psi(s,Y)^2 \dd s
        \bigg), \label{eq:1:generic_sv_model} \tag{SV1} \\
        Y_t
        &=
        Y_0 + \int_0^t K(t,s) \mu(s,Y) \dd s + \int_0^t K(t,s) \sigma^Y(s,Y)\, \dd W_s
        \label{eq:2:generic_sv_model}, \tag{SV2}
    \end{align}
    for $S_0 > 0$ and $Y_0 \in \R^m$. 
\end{definition}

\begin{remark}
We remark that Definition~\ref{def:generic_sv_model} is quite generic. 
        Suppose that strong existence holds for  \eqref{eq:2:generic_sv_model}. This allows us to start on any filtered probability space satisfying the usual conditions and supporting a Brownian motion, and obtain a continuous process~$Y$ satisfying \eqref{eq:2:generic_sv_model}. 
        On the other hand, if  weak existence holds for \eqref{eq:2:generic_sv_model}, then the filtered probability space and the Brownian motion are part of the solution. 
Note that uniqueness of the solution
        of~\eqref{eq:2:generic_sv_model} is not relevant.
\end{remark}

The following results provide tools for establishing~\eqref{eq:sup} for various stochastic volatility models, presented in Section~\ref{section:sup_rBergomi}. We start with a quite general theorem as our basic result. It bounds \eqref{eq:sup} in terms of the expected integrated squared variance under a new measure~$\P'$. 

\begin{theorem} \label{thm:int_of_sup}
    Assume the setting of Definition~\ref{def:generic_sv_model}; in particular, we are given continuous processes $S$ and $Y$ satisfying \eqref{eq:1:generic_sv_model} and \eqref{eq:2:generic_sv_model}, respectively. Let~$\tau$ be any stopping time such that~$S^\tau$ is a uniformly integrable martingale. 
    Then there exists a probability measure $\P'$ on $(\Omega, \F)$  and an $\R^{d}$-valued $\P'$-Brownian motion $W'$ such that $Y$ satisfies
        \begin{equation} \label{eq:int_of_sup}
        \begin{aligned}
            Y_t
            =
            Y_0 
            &+ \int_0^t K(t,s) \mu(s,Y) \dd s \\
            &+ \int_0^{t \land \tau} K(t,s) \sigma^Y(s,Y)\sigma^S(s,Y)^\top 
            \psi(s,Y)
            \dd s 
            \\
            &+ \int_0^t K(t,s) \sigma^Y(s,Y) \, \dd
            W'_s, \quad t \geq 0.
        \end{aligned}
    \end{equation}
Moreover, we have the estimate
    \begin{align} \label{eq:251207}
        \Ebig{\sup_{t \in [0,\tau]} S_t}
		&\leq
        S_0\frac{e}{e-1}\bigg(1 
        + 
        \mathbb{E}_{\P'}\bigg[
             \int_0^\tau \psi(s,Y)^2 \dd s
        \bigg]\bigg).
    \end{align}
\end{theorem}

\begin{proof}
    Assume w.l.o.g.\ that $S_0 = 1$ and  
    define the (stopped) share measures
    ${\dd\P'} = S_\infty^{\tau} \dd\mathbb{P}$. Applying Girsanov's theorem yields that
    \begin{align*}
		W'
		=
        W
        -
		\int_0^{.\land \tau} \sigma^S(s,Y)^\top\psi(s,Y) \dd s
	\end{align*}
    is an $\R^{d}$-valued Brownian motion under $\P'$, and hence \eqref{eq:int_of_sup}. For the second assertion we assume w.l.o.g.\ that 
    \begin{align*}
        \mathbb{E}_{\P'}\bigg[
             \int_0^\tau \psi(s,Y)^2 \dd s
        \bigg] < \infty,
    \end{align*}
    since~\eqref{eq:251207} would otherwise be trivially satisfied. Thus, the random variable $\int_0^\tau \psi(s,Y)\sigma^S(s,Y)\dd W'_s$ has zero expectation under $\P'$. Note that a simple rearrangement shows for $t \geq 0$
    \begin{align*}
        \int_0^t \psi(s,Y) \sigma^S(s,Y)\dd W_s
        &=
        \int_0^t \psi(s,Y) \sigma^S(s,Y)\dd W'_s
        +
        \int_0^{t \land \tau} \psi(s,Y)^2 \dd s .
    \end{align*}
    Combining this with an application of Lemma~\ref{lem:exp_sup_est} to the uniformly integrable martingale $S^\tau$ yields
    \begin{align*}
        \Ebig{\sup_{t \in [0,\tau]} S_t}
		&\leq
        c\bigg(1 
        + 
        \Ebig{S_\tau
            \bigg(
                \int_0^\tau \psi(s,Y) \sigma^S(s,Y)\dd W_s - \frac12 \int_0^\tau \psi(s,Y)^2 \dd s
            \bigg)
        }\bigg) \\
        &=
        c\bigg(1
        + 
        \mathbb{E}_{\P'}\bigg[
            \int_0^\tau \psi(s,Y) \sigma^S(s,Y)\dd W'_s + \frac12 \int_0^\tau \psi(s,Y)^2 \dd s
        \bigg]\bigg) \\
        &=
        c\bigg(1 
        + 
        \mathbb{E}_{\P'}\bigg[
            \frac12 \int_0^\tau \psi(s,Y)^2 \dd s
        \bigg]\bigg),
    \end{align*}
    where $c = e/(e-1)$. This concludes the proof.
\end{proof}

Now we combine Theorem~\ref{thm:int_of_sup} with a straightforward localization argument.

\begin{corollary} \label{thm:2:int_of_sup}
    Assume the setup of Definition~\ref{def:generic_sv_model} 
    and let $(\tau_n)_{n \in \N}$ be a localizing sequence for
    the local martingale~$S$. 
    Then, for any $T > 0$,
    \begin{align} 
        \Ebig{\sup_{t \in [0,T]} S_t}
		&\leq
        S_0\frac{e}{e-1}
        \bigg(
        1
        + 
        \lim_{n \uparrow \infty}
        \mathbb{E}_{\P'_n}\bigg[
             \int_0^{\tau_n \land T} \psi(s,Y)^2 \dd s
        \bigg] \label{cor:int_of_sup:eq1} 
        \bigg),
    \end{align}
    where $\dd \P'_n = S^{\tau_n}_T/S_0 \dd \P$. 
    For each $n \in \N$, the process $Y$ satisfies \eqref{eq:int_of_sup} with $\tau = \tau_n \land T$ and  $W' = W'_n$.
\end{corollary}

\begin{proof}
    This is an immediate consequence of Theorem~\ref{thm:int_of_sup}, since $S^{\tau_n \land T}$ is a uniformly integrable martingale for each $n \in \N$. An application of monotone convergence yields \eqref{cor:int_of_sup:eq1}.
\end{proof}

Theorem~\ref{thm:int_of_sup} and Corollary~\ref{thm:2:int_of_sup} give two ways to show~\eqref{eq:sup} for a given model. If we know a priori that~$S$ is a martingale, then
the deterministic stopping time~$T$ turns $S^T$ into a uniformly integrable martingale. Thus, Theorem~\ref{thm:int_of_sup} gives a concrete upper bound we can check; for an example of this approach, see the proof of Theorem~\ref{thm:sig}. On the other hand, if we do not know whether~$S$ is a martingale, then we first choose a localizing sequence, e.g., given by $\tau_n = \inf\{t \geq 0: |Y_t| \geq n\}$ for each $n \in \N$, and then try to bound the right hand side
of~\eqref{cor:int_of_sup:eq1} uniformly in~$n$. If this works, then~$S$ is a martingale with integrable supremum. Theorem~\ref{thm:rBergomi_multi} provides an example for this approach.  See also~\cite{Ru15} for a related test of the martingale property, using an appropriate version of F\"ollmer's measure.

\section{Applications}
\label{section:sup_rBergomi}

We apply the results of Section~\ref{se:int}, in particular Theorem~\ref{thm:int_of_sup} and Corollary~\ref{thm:2:int_of_sup}, to the validity of~\eqref{eq:sup} for three (classes of) stochastic volatility models. First we treat the mixed multi-factor rough Bergomi model from~\cite{lacombe2020asymptotics}, which contains the classical rough Bergomi model as a special case. We continue with a class of models whose coefficients satisfy certain growth assumptions, which are in particular satisfied by the rough Heston model. Lastly, we show integrability of the supremum for multi-factor signature volatility models driven by multiple Brownian motions, as treated in~\cite{jaber2025martingalepropertymomentexplosions}.

\subsection{Rough Bergomi model and multi-factor extensions}

We introduce the multi-factor extension of the rough Bergomi model from~\cite{lacombe2020asymptotics}. This definition extends the classical rough Bergomi model from~\cite{BaFrGa16} by incorporating multiple driving Brownian motions. To this end, we first define the fractional Riemann-Liouville kernel 
\begin{align} \label{def:kernel}
	K_{\alpha, \eta}(t,s) = 
		\eta \sqrt{2\alpha - 1} (t-s)^{\alpha - 1},  
  \quad  
  0 \leq s < t,
\end{align}
with $\alpha > 1/2$ and $\eta > 0$. Further, for $m \in \N$ consider a $d=(m+1)$-dimensional Brownian motion $W = (W^1,\dots,W^d)$ 
and define the vector-valued Riemann-Liouville process $Y$ by  
    \begin{align*}
    	Y_t
    	= 
        (Y^1_t, \dots, Y^m_t) 
    	= 
        \bigg(
            \int_0^t K_{\alpha, \eta} (t,s) \dd W^1_s, \dots, \int_0^t  K_{\alpha, \eta} (t,s) \dd W^m_s
        \bigg),
        \quad t \geq 0.
    \end{align*}

\begin{definition}[Mixed multi-factor rough Bergomi model \cite{lacombe2020asymptotics}] \label{model:rBergomi_multi}
Let $N \in \N$. 
    The mixed multi-factor rough Bergomi model is given in terms of the stock price and its instantaneous variance process for $t \geq 0$ by
    \begin{align*}
        S_t 
        &= 
        S_0 \exp\left(\int_0^t \psi(s,Y) \sigma^S \dd W_s - \frac12 \int_0^t \psi(s,Y)^2 \dd s\right), \\
        \psi(t, Y)
        &=
        \bigg(\xi_0(t) \sum_{i=1}^N \gamma_i \exp
        \Big(
            \frac{\nu_i^\top L_i}{\eta} Y_t - \frac{\nu_i^\top \Sigma_i \nu_i}2 t^{2\alpha-1}
        \Big)\bigg)^{\frac12},
    \end{align*}
    where $S_0 > 0$, $\xi_0: [0,\infty) \to (0,\infty)$ is some continuous and bounded function representing the initial forward variance curve, and $\sigma^S$ is a deterministic $d$-dimensional vector with $\|\sigma^S\|_2=1$. 
    Furthermore, $(\gamma_1, \dots, \gamma_N) \in (0,1]^N$ is such that $\sum_{i=1}^N \gamma_i = 1$, and the vector $\nu_i = (\nu_i^1, \dots, \nu_i^m)^\top$ satisfies $0 < \nu_i^1 < \dots < \nu_i^m$ for all $i \in \{1,\dots,N\}$. The matrix $L_i \in \R^{m \times m}$ is lower triangular and such that $\Sigma_i = L_i L_i^\top$ is positive definite for each $i \in \{1,\dots,N\}$.
\end{definition} 

The next theorem fully characterizes the martingale property and integrability of the supremum~\eqref{eq:sup} for the mixed multi-factor rough Bergomi model. The result regarding the martingale property follows the same spirit as~\cite{Ga19}. The characterization of~\eqref{eq:sup} is novel. Due to its significance in the literature we give a detailed remark at the end of this section on how the standard rough Bergomi model fits this definition and how the main theorem of this section applies in this situation.  

\begin{theorem} \label{thm:rBergomi_multi}
    Consider the mixed multi-factor rough Bergomi model from
    Definition~\ref{model:rBergomi_multi}, and denote by $\sigma^S_{(1:m)}$ and $\sigma^S_{(m+1)}$ the split of $\sigma^S$ into its first $m$ and last components, respectively. 
    Then the stock price $S$ is a martingale if and only if $\sigma^S_{(1:m)} L_i^\top \nu_i  \leq 0$ for all $i =1,\dots,N$. In this case, \eqref{eq:sup} holds for all $T > 0$.
\end{theorem}

\begin{proof}
Clearly, the mixed multi-factor rough Bergomi model fits Definition~\ref{def:generic_sv_model}. 
  We first show that $S$ is a martingale and \eqref{eq:sup} holds under the stated assumption.
    To this end, let~$(\tau_n)_{n \in \N}$ be a localizing sequence such that~$S^{\tau_n}$ is a bounded martingale. 
    Fix now $n \in \N$. By Corollary~\ref{thm:2:int_of_sup},  there exists a probability measure $\P'_n$ and a $\P'_n$-Brownian motion $W'_n$ such that
    \begin{align*}
        Y_t 
        = 
        \int_0^{t \land \tau_n \land T} K_{\alpha,\eta}(t,s) 
        \psi(s,Y) (\sigma^S_{(1:m)})^\top
        \dd s
        +
        \int_0^t K_{\alpha,\eta}(t,s) \dd W'_{n, 1:m,s}, \quad t \geq 0.
    \end{align*}
    For convenience, we write $M_n$ for the Gaussian process defined by $M_{n,t} = \int_0^t K_{\alpha,\eta}(t,s) \dd W'_{n, 1:m,s}$. 
    Using this representation of $Y$ and the assumption $\sigma^S_{(1:m)} L_i^\top \nu_i \leq 0$ for all $i =1,\dots,N$, 
    for any $t \in [0,T]$ one has 
     \begin{align*}
        \psi(t, Y)^2 
        &\leq 
        C_T \sum_{i=1}^N  
        \exp
        \bigg(
            \frac{\nu_i^\top L_i}{\eta}
            M_{n,t} 
            - 
               \frac{\nu_i^\top \Sigma_i \nu_i}2 t^{2\alpha-1}
        \bigg),
    \end{align*}   
    where $C_T > 0$ is some appropriate constant.     Thus one gets
    \begin{align*}
        &\mathbb{E}_{\P'_n}\bigg[
             \int_0^{\tau_n \land T} \psi(s,Y)^2 \dd s
        \bigg] \\
        &\hspace{50pt}\leq
        C_T \sum_{i=1}^N 
        \int_0^{T} \mathbb{E}_{\P'_n}\bigg[\exp
        \Big(
            \frac{\nu_i^\top L_i}{\eta}
            M_{n,s}
            -\frac{\nu_i^\top \Sigma_i \nu_i}2 s^{2\alpha-1}
        \Big)\bigg] \dd s
        = C_T N T.
    \end{align*}
 Hence $S$ is a martingale with integrable supremum by Corollary~\ref{thm:2:int_of_sup}.
 The converse implication is more involved and deferred to Appendix~\ref{app:sec:proof_thm:rBergomi_multi}. 
\end{proof}

\begin{remark}[Standard rough Bergomi model~\cite{BaFrGa16}]
    In case of the standard rough Bergomi model~\cite{BaFrGa16}, Definition~\ref{model:rBergomi_multi} and Theorem~\ref{thm:rBergomi_multi} simplify significantly. Indeed, choosing $m=N=\gamma_1=L_1=1$ and $\nu_1^1=\eta$ gives the instantaneous variance
    \begin{align*}
        \xi_0(t) \exp
        \bigg(
             \int_0^tK_{\alpha,\eta}(t,s)\dd W^1_s - \frac{\eta^2}2 t^{2\alpha-1}
        \bigg).
    \end{align*}
    Setting $\sigma^S \equiv (\rho, \; \sqrt{1-\rho^2})$ for $\rho \in [-1,1]$ yields a correlation $\rho$ between the Brownian motion driving the stock $S$ and $Y$, i.e., $\sigma^S W = \rho  W^1 + \sqrt{1-\rho^2}  W^2$. 
    In this case, the condition in Theorem~\ref{thm:rBergomi_multi} reduces to $\rho \leq 0$, which is thus a necessary and sufficient condition for the integrability of the supremum in the rough Bergomi model. Regarding the characterization of the martingale property, this exactly recovers the results in~\cite{Ga19}.
\end{remark}

\begin{remark}[General Gaussian Volterra models]
    The assumption that the variance is an exponential of~$Y$ is not essential to show that $S$ is a martingale with integrable supremum. Many papers, such as~\cite{BoJaMu25}, consider models with more general linking functions. In the uni-variate case, suppose we replace the linking function $(t,y) \mapsto \exp(y - \eta^2/2 t^{2\alpha-1})$ in the rough Bergomi model by a more general function $\zeta(t,y)$. As long as $\zeta$ is non-negative and increasing in the $y$-argument, and  the  first component of $\sigma^S$ is non-positive, it suffices to check that 
    \[
       \int_0^T\E{\zeta\Big(t, \int_0^tK_{\alpha,\eta}(t,s)dW_s\Big)} \dd t<\infty.
    \]
\end{remark}

\subsection{Models with coefficients of linear growth}

We apply the results of Section~\ref{se:int} to stochastic volatility models where the coefficient functions of Definition~\eqref{def:generic_sv_model} satisfy certain growth conditions.

\begin{assumption}[Linear growth] \label{ass:growth}
    Let $(S,Y)$ be a generic stochastic volatility model in the sense of Definition~\ref{def:generic_sv_model}. We assume that there exists an increasing function $D:[0,\infty) \to [0,\infty)$ such that
    \begin{align*}
        |\mu(t,y)| + |\psi(t,y)| + |\sigma^Y(t,y)| + |\sigma^Y(t,y) \sigma^S(t,y)^\top
        \psi(t,y)|
        &\leq D_t \left(1 + \sup_{s \leq t}|y_s|\right)
    \end{align*}
    holds for any $t \geq 0$ and $y \in C([0,\infty),\R^m)$
    and one of the following two conditions is satisfied:
    \begin{itemize}
        \item[(A)]
        The kernel $K$ satisfies $K(t,s) = 1$ for all $(t,s) \in \Delta$.
        \item[(B)]
        Assume that the functions $\mu, \sigma^S, \sigma^Y, \psi$ depend at time $t \geq 0$ on $Y$ only through the value $Y_t$ and in the growth assumptions from above $\sup_{s \leq t}|y_s|$ is replaced by $|y_t|$. Additionally assume there exists $\overline{K} \in L^2_{\mathrm{loc}}([0,\infty),\R)$ such that $K(t,s) = \overline K(t-s)$, i.e., $K$ depends on $s,t$ only through the difference $t-s$.
    \end{itemize}
\end{assumption}
Case (A) retains path dependency, but does not allow for singular kernels. Case~(B) allows for path dependency only through the kernel~$K$, which is important for applying the results of~\cite{AbLaPu19}. Using this assumption on generic stochastic volatility models, we can always conclude that~$S$ is a true martingale with integrable supremum by means of Corollary~\ref{thm:2:int_of_sup}.

\begin{theorem} \label{thm:pqintsup}
    Let $(S,Y)$ define a generic stochastic volatility model satisfying Assumption~\ref{ass:growth}. Then we have
    \begin{align*}
        \Ebig{\sup_{t \in [0,T]}S_t} < \infty
    \end{align*}
    for any $T > 0$. In particular, $S$ is a true martingale with integrable supremum.
\end{theorem}

\begin{proof}
    Fix $T > 0$ and define the stopping times $\tau_n = \inf\{t \geq 0: |Y_t| \geq n\}$ for each $n \in \N$. By Novikov's condition this yields a localizing sequence for~$S$.  
    Define next the stopped share measures $\dd \P_n' = S_T^{\tau_n}/S_0 \dd \P$ for each $n \in \N$. 
    Thanks to Corollary~\ref{thm:2:int_of_sup}, finiteness of the
    expression in~\eqref{cor:int_of_sup:eq1} implies the assertion. By the linear growth assumption, it suffices to show, for case~(A), that
    \begin{align*}
        \lim_{n \uparrow \infty}         \mathbb{E}_{\P'_n}\Big[\sup_{t \leq T} |Y_{t \land \tau_n}|^2\Big] < \infty,
    \end{align*}
    and, for case~(B), that
    \begin{align*}
        \lim_{n \uparrow \infty}         
         \int_0^T \mathbb{E}_{\P'_n}[|Y_t|^2 \1_{\{t \leq \tau_n\}}] \dd s  
         < \infty.
    \end{align*}
    In particular, without loss of generality, we may assume that $Y_0 = 0$. To make headway, note that for each $n \in \N$, thanks to Corollary~\ref{thm:2:int_of_sup}, we have the representation
    \begin{equation}
    \begin{aligned} \label{eq:lin_growth_1}
        Y_t
        =
 \int_0^{t \land \tau_n} K(t,s)\bar \mu(s,Y) \dd s
        &+ \int_{t \land \tau_n}^{t} K(t,s)\mu(s,Y) \dd s \\
        &+ \int_0^t K(t,s) \sigma^Y(s,Y) \, \dd W_{n,s}' 
    \end{aligned}
    \end{equation}
    for $t \in [0,T]$ and $\P_n'$-Brownian motions $W_n'$, where $\bar \mu = \mu + \sigma^Y (\sigma^S)^\top \psi $ is of linear growth by assumption.

    First assume that case~(A) of Assumption~\ref{ass:growth} is in force and recall
    that $K \equiv 1$ in this case. 
    Fix $n \in \N$ for the moment. Then the BDG inequality, together with Jensen's inequality, yields the existence of some $C_T > 0$, independent of $n$, such that
    \begin{align*}
        \mathbb{E}_{\P'_n}\Big[\sup_{s \leq T} |Y_{s \land \tau_n}|^2\Big]
        &\leq 
        C_T \bigg(1 + 
        \int_0^T 
        \mathbb{E}_{\P'_n}\big[|\bar\mu(t ,Y)|^2 + |\sigma^Y(t,Y)|^2\big] \dd t\bigg).
    \end{align*}
    Using the linear growth assumptions on $\bar \mu$ and $\sigma^Y$ yields the existence of another constant $C_T'$, independent of $n$, such that
    \begin{align*}
        \mathbb{E}_{\P'_n}\Big[\sup_{t \leq T} |Y_{t \land \tau_n}|^2\Big]
        &\leq
        C_T'\bigg(
       1 + 
        \int_0^T
        \mathbb{E}_{\P'_n}\Big[\sup_{s \leq t}|Y_{s \land \tau_n}|^2\Big] \dd t \bigg).
    \end{align*}
 Now Gr\"onwall's inequality    implies 
    \begin{align*}
    \mathbb{E}_{\P'_n}\Big[\sup_{t \leq T} |Y_{t \land \tau_n}|^2\Big]
        \leq 
        C_T' e^{C_T' T}, 
    \end{align*}
yielding case~(A). 

    Assume now that case~(B) is in force. We follow the argument from \cite[Proof of Lemma~3.1]{AbLaPu19} with slight alterations. Due to the assumption that all coefficient functions depend at time $s$ on the path $Y$ only through $Y_s$, we abuse notation and view them as function of $(s,Y_s)$. We fix again $n \in \N$. 
    We first claim that for $t \in [0,T]$,
    \begin{equation}
        \begin{aligned} \label{eq:1:growth_thm}
            |Y_t|^2 \1_{\{t < \tau_n\}}
            \leq 
            \bigg|
                &\int_0^t K(t,s) \bar \mu(s, Y_s \1_{\{s < \tau_n\}})
                \dd s \\
                &+ \int_0^t K(t,s) \sigma^Y(s, Y_s \1_{\{s < \tau_n\}}) \dd 
                W_{n,s}'
            \bigg|^2.
        \end{aligned}
    \end{equation}
    Indeed, on the event $\{t \geq \tau_n\}$, the left hand side is zero and the right hand side is nonnegative. 
    Moreover, on the event $\{t < \tau_n\}$, we have, thanks to It\^o isometry, 
    \begin{align*}
        \int_0^tK(t,s) \sigma^Y(s,Y_s) \dd W_{n,s}'
        =
        \int_0^{t}K(t,s) \sigma^Y(s,Y_s \1_{\{s < \tau_n\}}) \dd W_{n,s}'
    \end{align*}    
    for all $t \in \mathbb{Q} \cap [0, T]$, hence by continuity also for all $t \in [0, T]$. A similar statement holds for the finite-variation integral on the right-hand side
    of~\eqref{eq:1:growth_thm}.  Recalling \eqref{eq:lin_growth_1} thus
    yields~\eqref{eq:1:growth_thm}.
    
     Using Jensen's inequality on~\eqref{eq:1:growth_thm}, there exists some $C_T > 0$, independent of $n \in \N$, such that, for each $t \in [0,T]$,
    \begin{align*}
        |Y_t|^2 \1_{\{t < \tau_n\}}
        \leq
        C_T \bigg(
       1
        &+
        \int_0^t K(t,s)^2 |\bar \mu(s, Y_s \1_{\{s < \tau_n\}})|^2 \dd s \\
        &+ 
        \bigg|\int_0^t K(t,s) \sigma^Y(s, Y_s \1_{\{s < \tau_n\}}) \dd W_{n,s}'\bigg|^2 \bigg).
    \end{align*}
    Writing $f_n(t) = \mathbb{E}_{\P_n'}[|Y_t|^2 \1_{\{t < \tau_n\}}]$  for each $t \in [0, T]$ and applying the BDG inequality to the $\P'_n$-local martingale $(M_r)_{r \in [0,t]}$ with 
    \[
      M_r = \int_0^r K(t,s) \sigma^Y(s, Y_s \1_{\{s < \tau_n\}}) \dd W_{n,s}'
    \]
    yields, using the linear growth assumption, a constant $C_T'$, independent of $n \in \N$, such that
    \begin{align*}
        f_n(t)
        &\leq
        C_T \Big(
        1
        +
        \int_0^t K(t,s)^2 \mathbb{E}_{\P'_n}[|\bar \mu(s, Y_s \1_{\{s < \tau_n\}})|^2] \dd s
        +
        \mathbb{E}_{\P'_n}\Big[\sup_{r \in [0,t]} |M_r|^2\Big] \Big)\\
        &\leq
        C_T' \Big(
            1 + \int_0^T \overline K(s)^2 \dd s +  \int_0^t K(t,s)^2 f_n(s) \dd s
        \Big).
    \end{align*}
    We are now exactly in the situation $f_n \leq c' + c' K^2 * f_n$ for some $c' > 0$ on $[0,T]$. Following the same arguments as in \cite[Proof of Lemma~3.1]{AbLaPu19}, we have $\sup_{t \in [0,T], n \in \N}f_n(t) \leq c$ for some $c > 0$, concluding the proof.
\end{proof}

\begin{remark}
    The last theorem in particular applies to affine Volterra processes in the sense of \cite{AbLaPu19} and hence also to the well-known rough Heston model. Note however that the existence of higher order moments in the rough Heston model is known, see \cite{GeGePi19}. In this case, Theorem~\ref{thm:pqintsup} gives an alternative proof of~\eqref{eq:sup}.
\end{remark}

\subsection{Signature volatility models}

Signature volatility models have recently attracted considerable attention as a flexible, nonparametric framework for modeling volatility in financial markets. Drawing on ideas from rough path theory and signature methods, they offer a way to include path-dependency in a tractable way that classical parametric models can not. Notable work includes~\cite{cuchiero2023spvix}, which develops a signature-based model to jointly calibrate SPX and VIX options across maturities, and \cite{abijaber2024signature}, which uses signature based models for pricing and hedging derivatives in models based on signature of Brownian motion. In particular, for Brownian motion based signature volatility models, \cite{jaber2025martingalepropertymomentexplosions} characterizes the martingale property and provides conditions under which higher order moments exist. 

Before elaborating these results we give a minimalist introduction to signature. For a more thorough introduction, see for example \cite{abijaber2024signature, cuchiero2023spvix, cuchiero2022theocalib}. Since it suffices for our purposes, we only introduce truncated signatures and fix a truncation level $N \in \N$ and some dimension $d \in \N$ throughout this subsection. For $n \in \N$, we denote the $n$-th tensor power of $\R^d$ by $(\R^d)^{\otimes n}$ with $(\R^d)^{\otimes 0} = \R$. Using this, we define the truncated tensor algebra of order $N$ as
\begin{align*}
    T^N(\R^d)
    &=
    \{
    \mathbf{a} =(a_0,\dots,a_N): a_n \in (\R^{d})^{\otimes n}, n =0,\dots,N
    \}.
\end{align*}
For each $\mathbf{a}, \mathbf{b} \in T^N(\R^d)$ and $\lambda \in \R$ we define scalar multiplication and addition 
\begin{align*}
    \lambda \cdot \mathbf{a} = (\lambda a_0, \dots,\lambda a_N), \quad 
    \mathbf{a} + \mathbf{b} = (a_0 + b_0, \dots,a_N + b_N)
\end{align*}
and the tensor product
\begin{align*}
    \mathbf{a} \otimes \mathbf{b}
    = 
    (c_0, \dots, c_N), \quad 
    c_n = \sum_{k=0}^n a_k \otimes b_{n-k}.
\end{align*}
For $n \in \{0, \ldots, N\}$, denote the set of words of length~$n$ over the alphabet $\{1,\dots,d\}$ by $\mathcal{W}_n$. 
Note that the set of words of length zero $\mathcal{W}_0$, contains only the empty word $\emptyset$. For $n \in \{0, \ldots, N\}$ and a word $w = (i_1,\dots,i_n)$ we set $e_w = e_{i_1}\otimes \dots \otimes e_{i_n}$ such that each element $\mathbf{a} \in T^N(\R^d)$ can be written as 
\begin{align*}
    \sum_{n=0}^N\sum_{w \in \mathcal{W}_n} a_w e_w, \quad a_w \in \R.
\end{align*}
For $\mathbf{a} =  \sum_{n=0}^N\sum_{w \in \mathcal{W}_n} a_w e_w \in T^N(\R^d)$ and $\mathbf{b} =  \sum_{n=0}^N\sum_{w \in \mathcal{W}_n} b_w e_w\in T^N(\R^d)$ we set 
\begin{align*}
    \ran{\mathbf{a}}{\mathbf{b}}
    =
    \sum_{n=0}^N\sum_{w \in \mathcal{W}_n} a_w b_w.
\end{align*}
Observe in particular that $a_w = \ran{e_w}{\mathbf{a}}$ for any word $w$. 

Let $X$ now be a continuous semimartingale. Then its truncated signature is recursively defined as the $T^N(\R^d)$-valued process $\mathbb{X}$ satisfying
\begin{align*}
    \ran{e_\emptyset}{\mathbb{X}_t} = 1, \quad
    \ran{e_w}{\mathbb{X}_t} 
    =
    \int_0^t \ran{e_{w'_n}}{\mathbb{X}_t} \circ \dd X^{i_n}_t, \quad w = (i_1,\dots,i_n) \in \mathcal{W}_n
\end{align*}
for $n \in \{ 0,\ldots, N\}$, where $\circ$ denotes Stratonovich integration, and $w'_n$ is defined as $(i_1,\dots,i_{n-1})$ for the word $w = (i_1,\dots,i_n)$. Using the definition of the tensor product we see that the truncated signature of~$X$ can equivalently be defined as the unique solution to 
\begin{align} \label{eq:sig_SDE}
    \dd \mathbb{X}_t = \mathbb{X}_t \otimes \circ \dd X_t, \quad \mathbb{X}_0 = e_\emptyset.
\end{align} 

We introduce the multi-factor signature volatility model from \cite[Section 3.3]{jaber2025martingalepropertymomentexplosions}. 
Consider a $d$-dimensional Brownian motion $W = (W^1,\dots,W^d)$ on a filtered probability space satisfying the usual conditions. 
The signature volatility model is defined as 
\begin{equation} \label{eq:sig_vol_model}
    \begin{aligned} 
        \frac{\dd S_t}{S_t} 
        &= 
        \xi_t \dd W^{d}_t, \\
        \xi_t
        &= 
        \ran{\bfsig}{\X_t},
    \end{aligned}
\end{equation}
where $\X$ is the truncated signature of the time extended Brownian motion $(t,W)$ and $\bfsig \in T^N(\R^{d+1})$. The question of integrability of the supremum~\eqref{eq:sup} has been answered in \cite[Theorem~3.2]{jaber2025martingalepropertymomentexplosions} for the case $m=1$ via existence of higher moments. The existence of higher moments for the model~\eqref{eq:sig_vol_model} driven by multiple Brownian motions is not addressed in \cite{jaber2025martingalepropertymomentexplosions}. Hence the question of~\eqref{eq:sup} is still open, and the next result directly characterizes \eqref{eq:sup} in this case.  

\begin{theorem} \label{thm:sig}
    Consider the multi-factor signature volatility model~\eqref{eq:sig_vol_model} with $N \geq 2$ and define $\sigma' = \ran{e_{d+1}^{\otimes N}}{\bfsig}$, which we assume to be non-zero. Then~\eqref{eq:sup} holds if and only if $N$ is odd and $\sigma' < 0$. 
\end{theorem}

\begin{proof}
    Under our assumptions, the process $S$ is a martingale if and only if $N$ is odd and $\sigma' < 0$, by 
    \cite[Theorem~3.5]{jaber2025martingalepropertymomentexplosions}. 
    Hence it suffices to show \eqref{eq:sup} if $N$ is odd and $\sigma' < 0$.

    Denote $m = \dim T^N(\R^{d+1})$.  Then $T^N(\R^{d+1})$ is isomorphic to $\R^m$ as a vector space.
    Define $\psi(t,\mathbf{y}) = \ran{\bfsig}{\mathbf{y}_t}$ for $t \geq 0$ and $\mathbf{y} \in C([0,\infty),\R^m)$, and $K = 1$. Noting that $\X$ satisfies the linear SDE \eqref{eq:sig_SDE}, with $X$ replaced by $(t,W)$, the signature volatility model fits into Definition~\ref{def:generic_sv_model}.

    The dynamics of $\X$, corresponding to \eqref{eq:int_of_sup}, under the measure $\P'$ from Theorem~\ref{thm:int_of_sup} are now given by
    \begin{equation} \label{eq:1:sig_SDE}
        \begin{aligned} 
            \dd \X_t
            = 
            \big(\X_t \otimes e_1
            +
            \ran{\bfsig}{\X_t} \X_t \otimes e_{d+1}\big) \dd t
            +
            \sum_{i=2}^{d+1} \X_t \otimes e_i \circ \dd (W^{i-1})'_t.
        \end{aligned}
    \end{equation} 
    Consider the signature SDE
    \begin{align*} 
        \dd Y_t = \ran{\bfsig}{\mathbb{\widehat Y}_t} \dd t + \dd (W^d)'_t, \quad Y_0 = 0, \qquad t \in [0,T],
    \end{align*}
    where $\mathbb{\widehat Y}_t$ is the signature of the process $(t,W_{1:d-1}',Y)$. Note that the signature SDE \eqref{eq:sig_SDE} for $\mathbb{\widehat Y}_t$ agrees with \eqref{eq:1:sig_SDE}. Since the vector fields therein are locally Lipschitz, the solution is unique up to a time of explosion $\tau_\infty$. In \cite[Section~6]{abijaber2024signature} it is shown that $\P[\tau_\infty < T] = 0$ and thus, since $\P' \ll \P$ by construction, also $\P'[\tau_\infty < T] = 0$. Hence, we have the equality
    \begin{align} \label{eq:4:sig_SDE}
        \psi(t, \X_t)
        =
        \ran{\bfsig}{\mathbb{\widehat Y}_t}, \quad  t \in [0,T].
    \end{align}
    By \cite[Lemma~5.2 and Section~6]{jaber2025martingalepropertymomentexplosions}  there exists a constant $C > 0$ such that 
    \begin{align} \label{eq:2:sig_SDE}
        \mathbb{E}_{\P'}\Big[\int_0^{T}
            \ran{\bfsig}{\mathbb{\widehat Y}_t}^2 \dd t
        \Big]
        \leq
        C\Big(1 + 
        \mathbb{E}_{\P'}\Big[\int_0^{T}
           Y_t^{2N} \dd t
        \Big]\Big)
    \end{align}
    and 
    \begin{align} \label{eq:3:sig_SDE}
        \mathbb{E}_{\P'}[{Y^{2N}_{t}}] \leq C, \quad t \in [0,T].
    \end{align}
    Combining \eqref{eq:4:sig_SDE}, \eqref{eq:2:sig_SDE}, and \eqref{eq:3:sig_SDE} thus yields  
    \begin{align*}
        \mathbb{E}_{\P'}\Big[ 
            \int_0^{T}
            \psi(t,\X)^2 \dd t
        \Big]
        \leq 
        C \Big(
        1 +
            \int_0^{T}
            \mathbb{E}_{\P'}[Y_t^{2N}] \dd t
        \Big) < \infty.
    \end{align*}
    An application of Theorem~\ref{thm:int_of_sup} then concludes the proof.
\end{proof}

\section{Martingales with non-integrable supremum}\label{se:sup}

\subsection{General remarks}

By the Burkholder-Davis-Gundy inequalities, the set of local martingales
with integrable supremum coincides with~$\mathcal{H}^1$, the space
of local martingales~$(M_t)_{t\geq0}$ for which
\[
  \| M \|_1 = \mathbb{E}\big[[M]_\infty^{1/2}\big]
\]
is finite.
By dominated convergence, ``local'' can be dropped from this statement.
Clearly,
every element of~$\mathcal{H}^1$
is uniformly integrable.
For details and further interesting properties of~$\mathcal{H}^1$, such as its duality
to the space of BMO martingales, we refer to Section IV.4 in~\cite{Pr04}.

A family of \emph{discrete time} examples of martingales with \emph{non}-integrable
supremum is provided by Theorem~2 in~\cite{BlDu63}. In this section,
we give three different approaches that yield 
examples in continuous time.
Further examples can be found in~\cite{Lo16}.
First, recall that
any continuous nonnegative martingale~$(X_t)_{t \geq 0}$ with $X_0=1$
satisfies the reverse $L^1$ inequality
\begin{equation}\label{eq:gundy}
  \mathbb{E}\Big[\sup_{t \geq 0} X_t\Big] \geq 1 + \mathbb{E}[X_\infty (\log(X_\infty))^+].
\end{equation}
For a proof, we refer to p.~149 in~\cite{Du84};
see also~\cite{Gu69,Gu80} for earlier versions of similar
inequalities.
For the mere existence of martingales with non-integrable supremum, the following
observation suffices: For $T\in(0,\infty]$, suppose that a probability space with filtration $(\mathcal{F}_t)_{0\leq t\leq T}$ supports a nonnegative $\mathcal{F}_T$-measurable random variable~$X_T \in L^1$ which is
not in the class $L\log L$, defined at the end of the introduction. 
If  $X_t=\mathbb{E}[X_T|\mathcal{F}_t]$ is
continuous, then it is
 a uniformly integrable martingale which is not in $\mathcal{H}^1$.

\subsection{A construction by Dubins and Gilat}

The second construction is due to~\cite{DuGi78} and appears also 
in~\cite[Example~1]{Lo16}. In our exposition,
we add some details.
In contrast to the above argument, the process will be explicitly
defined.
Unlike the rest of the paper, we drop the assumption of non-negativity, and aim at 
finding a real-valued martingale $(X_t)_{0\leq t\leq 1}$ with $\mathbb{E}[\sup_{t\in[0,1]}X_t]=\infty$.
For convenience, the following assumption is in force throughout this subsection.
\begin{assumption}\label{ass:F}
Let~$F$ be the cumulative distribution function of some distribution on the real line with finite
first moment, with full support and a continuous density.
\end{assumption}

The function
\begin{equation}\label{eq:def H}
  H_F(t)= \frac{1}{1-t} \int_t^1 F^{-1}(s)\dd s,\quad 0< t<1,
\end{equation}
is known as the Hardy--Littlewood  maximal function of the distribution.
For the significance of this function concerning analytic and martingale
inequalities, we refer to~\cite{Gi86} and the references therein.
Note that, in financial terms, the Hardy--Littlewood maximal function is
the tail value at risk, or expected shortfall, of the distribution defined by~$F$.
In economics, the closely related function $t\mapsto\int_0^t F^{-1}(s)\dd s/ \int_0^1 F^{-1}(s)\dd s$ is known as the Lorenz curve~\cite{KaRa15}. 
It has been shown in~\cite{BlDu63} that the distribution
of~$H_F$, viewed as a random variable on the probability space
\begin{equation}\label{eq:ps}
 \big((0,1), \mathcal{B}(0,1), \mathrm{Leb} \big),
\end{equation}
is an upper bound for the set
\begin{multline*}
\{  \nu : \nu \text{ is the distribution of the supremum of a martingale closed} \\ 
\text{by a random variable with cumulative distribution function }F\}
\end{multline*}
with respect to stochastic order. In~\cite{DuGi78} it was shown that the distribution
of~$H_F$ is also a member of this set.
This is achieved by defining a martingale~$X$ on the probability space~\eqref{eq:ps}
such that its closing element has cumulative distribution function $F$ and the distribution of the supremum of~$X$ equals the distribution of $H_F$. If $H_F\notin L^1$, then~$X$ has the desired properties. See Proposition~\ref{prop:F} for details.

\begin{remark}\label{rem:diesundas}
  \begin{itemize}
    \item[(i)] It is clear that only the upper tail of the supremum of a martingale can cause 
    non-integrability, since
    \[
      \Big|\!\sup_{t\in[0,1]}X_t\Big|\ \mathbf{1}_{\{\sup_{t\in[0,1]}X_t<0\}}
      \leq |X_0|.
    \]
    \item[(ii)] Integrability of the Hardy--Littlewood maximal function~\eqref{eq:def H} can
    only fail because of a blowup as $t$ tends to $1$. Indeed,
    $\lim_{t\downarrow0}H_F(t)=\int_0^1 |F^{-1}(s)|\dd s$
    is finite, by existence of the first moment of~$F$. Thus, for integrability criteria 
    concerning~$H_F$, it does not matter
    if $|F^{-1}(s)|$ is used in the definition~\eqref{eq:def H}
    instead of $F^{-1}(s)$, as is the case in~\cite{St69}.
    \item[(iii)] Concerning notation: The nondecreasing function~$f$ from~\cite{DuGi78} is our~$F^{-1}$. On the other
    hand, $f^*$ from~\cite{Gi86} is assumed to be nonincreasing, and equals $F^{-1}(1-\cdot)$
    in our notation. As a consequence, when reading~\cite{DuGi78,Gi86} it is important
    to note that the Hardy--Littlewood maximal functions~$H$ from~\cite{DuGi78}
    (which is our~$H_F$) and~$F^*$
    from~\cite{Gi86} are not equal, but related by $F^*(1-t)=H(t)$, where $0<t<1$.
  \end{itemize}
\end{remark}

\begin{lemma}\label{le:stein}
If 
\begin{equation}\label{eq:x log x}
  \int_{\mathbb R} |x| (\log|x|)^+ F(\dd x) = \infty,
\end{equation}
then $H_F \notin L^1$.
\end{lemma}
\begin{proof}
   This follows from a classical theorem due to Stein~\cite{St69}, which asserts that
   $H_F \notin L^1$ is equivalent to $\int_{0}^1|F^{-1}(s)| (\log|F^{-1}(s)|)^+ \dd s= \infty$, by substitution.
\end{proof}
Examples of cumulative distribution functions with finite first moment, but satisfying \eqref{eq:x log x}, can be easily found; e.g., assume that the density satisfies
\[
  F'(x) \leq \frac{c_1}{x^2 (\log|x|)^\alpha}, \quad |x| \ \ \text{large,}
\]
and
\[
  F'(x) \geq \frac{c_2}{(x \log x)^2}, \quad x>0\ \ \text{large,}
\]
for some constants $\alpha \in (1,2]$ and $c_1,c_2>0$.

\begin{proposition}\label{prop:F}
  Recall Assumption~\ref{ass:F}, and suppose that~$F$ satisfies~\eqref{eq:x log x}.
  Then the process $(X_{t})_{0\leq t\leq 1}$  defined by
  \[
    X_t(s) =
    \begin{cases}
      H_F(t) & \quad \text{for}\quad t\leq s, \\
      F^{-1}(s) & \quad \text{for}\quad t>s
    \end{cases}
  \]
  is a martingale on the probability space~\eqref{eq:ps}, with respect to its own filtration,
  and its supremum is not integrable.
\end{proposition}
\begin{proof}
  According to the proof of~\cite[Lemma~2]{DuGi78}, the process $(X_{t})_{0\leq t\leq 1}$
  is a martingale. As this is not shown in~\cite{DuGi78}, we give a proof for the reader's convenience. Let~$U(s)=s$, so that~$U$ is a standard uniformly distributed random variable
  on the probability space~\eqref{eq:ps}. It is easy to see that
  $\mathcal{F}_t = \sigma(U \mathbf{1}_{\{U < t\}})$ is the filtration generated by~$X$.
  Define $\tilde{X}_t=\mathbb{E}[F^{-1}(U)|\mathcal{F}_t]$ for $t\in[0,1]$.
  It is readily verified that  we have $\tilde{X}_t = F^{-1}(U)$ on the event $\{U < t\} \in \mathcal{F}_t$, while on the event $\{U \geq t\} \in \mathcal{F}_t$ we have $\tilde{X}_t = H_F(t)$.
  Thus, $X=\tilde{X}$. Finally, we have
  \[
    \sup_{0\leq t\leq 1}X_t(s) = H_F(s), \quad 0<s<1.
  \]
  By Lemma~\ref{le:stein}, this is not integrable with respect to Lebesgue measure on $(0,1)$.
\end{proof}

\subsection{A new approach}

Our third construction starts with an arbitrary nonnegative local martingale~$M$
that is not a uniformly integrable martingale. In particular, $M$ may be any nonnegative
strict local martingale. See, e.g., Exercise 3.3.36 in~\cite{KaSh91} for a standard example.
We present the arguments for $T=\infty$, but they trivially apply to finite~$T$ as well.
Observe that\begin{equation*}
\mathbb{E}\Big[\sup_{t\geq 0}M_t\Big]=\int_0^\infty\mathbb{P}\Big[\sup_{t\geq 0}M_t>u\Big]\dd u
\end{equation*}
and
\begin{equation*}
\sum_{n=1}^\infty\mathbb{P}\Big[\sup_{t\geq 0}M_t>n\Big]
\leq\int_0^\infty\mathbb{P}\Big[\sup_{t\geq 0}M_t>u\Big]\dd u
\leq 1+\sum_{n=1}^\infty\mathbb{P}\Big[\sup_{t\geq 0}M_t>n\Big].
\end{equation*}
As mentioned above, we have $M\in\mathcal{H}^{1}$ if and only if $\mathbb{E}[\sup_{t\geq 0}M_t]<\infty$, and so it follows that $M\in\mathcal{H}^{1}$ if and only if
\begin{equation}\label{eq:nat cond}
\sum_{n=1}^\infty\mathbb{P}\Big[\sup_{t\geq 0}M_t>n\Big]<\infty.
\end{equation}
This condition plays a key role in our construction.

\begin{theorem}\label{thm:M}
  Let $(M_t)_{t\geq 0}$ 
  be a nonnegative local martingale that is not a uniformly integrable martingale. Then there is an extended filtered probability space with a stopping
  time~$\tau$ such that 
  \begin{itemize}
    \item[(i)] The lift of~$M$ to the extended probability space, again denoted by~$M$, has the same properties.
    \item[(ii)] The stopped process $M^\tau$ is a
    uniformly integrable martingale which
    is not in~$\mathcal{H}^1$, i.e., $\mathbb{E}[\sup_{t\geq0}M^\tau_t]=\infty$.
  \end{itemize}
\end{theorem}
\begin{proof}
  We assume that $T=\infty$; the construction
  clearly works for finite~$T$, too.
Define the nondecreasing sequence $(c_n)_{n\in\N}$, by setting
\begin{equation}
\label{eqSec2:1}
c_n= \log \bigg(e+\sum_{k=1}^n\mathbb{P}\Big[\sup_{t\geq 0}M_t >k\Big]\bigg) >1,	
\end{equation}
for each $n\in\N$. Since $M\notin\mathcal{H}^1$, it follows that $\lim_{n\uparrow\infty}c_n=\infty$. 
Without loss of generality, the original probability space,
with filtration $(\mathcal{F}_t)_{t\geq 0}$, accommodates an $\mathbb{N}$-valued 
and  $\mathcal{F}_0$-measurable random variable~$Y$ that is independent of $M$, and whose distribution satisfies $\mathbb{P}[Y>n]=1/c_n$, for each $n\in\N$. Let
\begin{equation}
\label{eqSec2:2}
\sigma=\inf\{t\geq 0\,|\,M_t>Y\}	
\end{equation}
denote the first time~$M$ exceeds~$Y$. It follows that 
\begin{equation*}
\mathbb{P}\Big[\sup_{t\geq 0} M^\sigma_t>n\Big]
\geq\mathbb{P}\Big[\sup_{t\geq 0}M_t>n\Big]\mathbb{P}[Y>n]
=\frac{1}{c_n}\mathbb{P}\Big[\sup_{t\geq 0}M_t>n\Big],
\end{equation*}
for each $n\in\N$. Consequently,
\begin{equation*}
\sum_{n=1}^\infty\mathbb{P}\Big[\sup_{t\geq 0}M^\sigma_t>n\Big]
\geq\lim_{m\uparrow\infty}\frac{1}{c_m}\sum_{n=1}^m\mathbb{P}\Big[\sup_{t\geq 0}M_t>n\Big]
=\lim_{m\uparrow\infty}\frac{e^{c_m}-e}{c_m}=\infty,
\end{equation*}
since $(c_n)_{n\in\N}$ is non-decreasing and $\lim_{n\uparrow\infty}c_n=\infty$. This implies that $M^\sigma\notin\mathcal{H}^1$. On the other hand, the almost sure limit $M^\sigma_\infty= M^\sigma_{\infty-}\in [0,\infty)$ exists, since $M^\sigma$ is a nonnegative local martingale, and hence also a nonnegative supermartingale. Moreover,
\begin{align*}
\mathbb{E}[M^\sigma_\infty]&=\sum_{n=1}^\infty\mathbb{E}[M^\sigma_\infty\,|\,Y=n]\mathbb{P}[Y=n]\\
&=\sum_{n=1}^\infty\mathbb{E}[M^{\tau_n}_\infty\,|\,Y=n]\mathbb{P}[Y=n]\\
&=\sum_{n=1}^\infty\mathbb{E}[M^{\tau_n}_\infty]\mathbb{P}[Y=n]
=\sum_{n=1}^\infty\mathbb{E}[M_0]\mathbb{P}[Y=n]
=\mathbb{E}[M^\sigma_0],
\end{align*}
where 
\begin{equation*}
\tau_n=\inf\{t\geq 0\,|\,M_t>n\}, \quad n \in \N.
\end{equation*}
Here, the penultimate equality follows from the fact that $M^{\tau_n}$ is
uniformly integrable, for each $n\in\N$. Consequently, $M^\sigma$ is uniformly integrable.
\end{proof}

Recall from the beginning of this section that the existence of martingales not in~$\mathcal{H}^1$
is an immediate consequence of~\eqref{eq:gundy}. If~$M$ is continuous, then the process~$M^\sigma$
from the proof of Theorem~\ref{thm:M} is of the kind mentioned there, i.e., 
$M_\infty^\sigma \notin L \log L$. Besides accommodating c\`adl\`ag processes~$M$,
we make the following pedagogical remark: Theorem~\ref{thm:M} and its proof require
only material from an introductory course on stochastic calculus, and not the reverse
$L^1$ inequality. Moreover, the natural way to prove $M^\sigma \notin \mathcal{H}^{1}$
is to show that~\eqref{eq:nat cond} fails, and not to verify $M_\infty^\sigma \notin L \log L$. 

\begin{example}
Consider a nonnegative local martingale~$M$ 
that belongs to Class~($\mathcal{C}_0$), according to the terminology of~\cite{NiYo06}, and suppose that $M_0=1$. In that case, $M$ is a strictly positive local martingale without any positive jumps, for which $M_\infty= M_{\infty-}=0$. The construction from Theorem~\ref{thm:M} is then applicable, since $\mathbb{E}[M_\infty]=0<1=\mathbb{E}[M_0]$ implies that~$M$ 
is not uniformly integrable. Moreover, an application of Doob's maximal identity~\cite[Lemma~2.1]{NiYo06} provides the following concrete representation for the non-decreasing sequence $(c_n)_{n\in\N}$, defined by~\eqref{eqSec2:1}:
\begin{equation*}
c_n=\log\biggl(e+\sum_{k=1}^n\frac{1}{k}\biggr), \quad n\in\N.
\end{equation*}
As the harmonic numbers tend to infinity,
we have $\lim_{n\uparrow\infty}c_n=\infty$, which is the crucial ingredient for showing that $M^\sigma\notin\mathcal{H}^1$, where the stopping
time~$\sigma$ is given by~\eqref{eqSec2:2}.
\end{example}

\appendix

\section{Auxiliary results for SVEs}
\label{app:SVIE_results}

This appendix provides auxiliary results needed in Appendix~\ref{app:sec:proof_thm:rBergomi_multi} for the proof of Theorem~\ref{thm:rBergomi_multi}.  Theorem~\ref{thm:VIE_comparison} establishes a comparison and uniqueness result for deterministic Volterra equations and  Theorem~\ref{thm:svie} shows existence and uniqueness for stochastic Volterra equations. 
We first recall the following definition from~\cite{gripenberg_londen_staffans_1990}. 

\begin{definition}[Volterra kernel of continuous type]  \label{D:241126}
    A measurable function $\kappa: \Delta \to [0, \infty)$ is a Volterra kernel of continuous type\footnote{By abuse of notation, for $t \in [0, \infty)$ we extend $\kappa(t,.)$ from $[0, t)$ to $[0, \infty)$ and set $\kappa(t,s) = 0$ for $s \geq t$.} if for each $t \in [0,\infty)$ we have $\kappa(t,.) \in L^1([0,\infty))$  and the map $[0,\infty) \ni t \mapsto \kappa(t,.)$ is continuous.
\end{definition}
To apply the results of~\cite{zhang2008stochastic} in the proof of the upcoming existence result for stochastic Volterra equations (SVEs), we establish two facts for Volterra kernels of continuous type. 
\begin{lemma} \label{app:lem1}
    Let $\kappa$ be a Volterra kernel of continuous type. Then for each bounded measurable function $g: \R \to \R$, the map 
    \begin{align*}
        t \mapsto \int_0^t \kappa(t,s) g(s) \dd s
    \end{align*}
    is continuous and for each $T > 0$ we have
    \begin{align*}
        \limsup_{\varepsilon \downarrow 0} 
        \Big\| 
        \int_\cdot^{\cdot+\varepsilon} \kappa(\cdot+\varepsilon, s) \dd s 
        \Big\|_{L^\infty([0,T])} 
        = 0.
    \end{align*}
\end{lemma}
\begin{proof}
    The first assertion is a special case of \cite[Theorem~9.5.3]{gripenberg_londen_staffans_1990},
    but we still prove it, as the following steps are required to prove the second assertion.
    Let $g: \R \to \R$ be bounded and measurable and let $h \geq 0$. Fix $t, u \geq 0$ and first assume $u > t$. Then 
    \begin{align*}
        &\Big| 
        \int_0^u \kappa(u+h,s) g(s) \dd s - \int_0^t \kappa(t+h,s) g(s) \dd s
        \Big| \\
        &\hspace{70pt}\leq 
        \|g\|_\infty \Big(
        \int_t^u \kappa(u+h,s)\dd s + \int_0^t |\kappa(u+h,s) - \kappa(t+h,s)|\dd s
        \Big).
    \end{align*}
    For $u \downarrow t$, the first integral on the right hand side vanishes via an application of dominated convergence
    and the second due to the assumed $L^1$-continuity. The case $u < t$ is similar. Hence, taking $h = 0$, the first claim is shown.
    
    For the second claim, fix $\varepsilon > 0$ and note that 
    \begin{align*}
        \int_t^{t+\varepsilon} \kappa(t+\varepsilon, s) \dd s
        = 
        \int_0^{t+\varepsilon} \kappa(t+\varepsilon, s) \dd s
        -
        \int_0^{t} \kappa(t+\varepsilon, s) \dd s
    \end{align*}
    is a continuous function of~$t$, since it is the difference of two continuous functions by the first part of the proof by taking $g = 1$ and $h = 0$ and $h = \varepsilon$, respectively. Hence, the supremum over the compact interval $[0,T]$ is attained by some maximizing $\hat t_\varepsilon \in [0,T]$. Hence, we have
    \begin{align*}
        \Big\| 
        \int_\cdot^{\cdot+\varepsilon} \kappa(\cdot+\varepsilon, s) \dd s 
        \Big\|_{L^\infty([0,T])} 
        = 
        \int_{\hat t_\varepsilon}^{\hat t_\varepsilon+\varepsilon} \kappa(\hat t_\varepsilon+\varepsilon, s) \dd s.
    \end{align*}
    Since the continuous maps
    \begin{align*}
        t \mapsto \int_0^t \kappa(t,s) \dd s \in \R
        \quad \text{and} \quad
        t \mapsto \kappa(t, \cdot) \in L^1([0,\infty))
    \end{align*}
    are uniformly continuous on the compact set $[0,T+1]$, there exists for every $\eta > 0$  some $\delta > 0$ such that $|u-v| < \delta$ implies 
    \begin{align*}
        \Big| \int_0^u \kappa(u,s) \dd s - \int_0^v \kappa(v,s) \dd s \Big| < \eta 
        \quad \text{and} \quad
        \int_0^\infty \big| \kappa(u,s) - \kappa(v,s) \big| \dd s< \eta. 
    \end{align*}
    Let $\eta > 0$ be arbitrary, choose $\delta$ accordingly, and fix an arbitrary $\varepsilon < \delta$. Then we have 
    \begin{align*}
        \Big|
        \int_{\hat t_{\varepsilon}}^{\hat t_{\varepsilon} +\varepsilon }\kappa({\hat t_{\varepsilon}}+\varepsilon, s) \dd s 
        \Big| 
        &= \Big|
        \int_0^{\hat t_{\varepsilon} +\varepsilon }\kappa(\hat t_{\varepsilon}+\varepsilon, s) \dd s - \int_0^{\hat t_{\varepsilon}}\kappa(\hat t_{\varepsilon}+\varepsilon, s) \dd s
        \Big| \\
        &\leq \Big|
        \int_0^{\hat t_{\varepsilon} +\varepsilon }\kappa(\hat t_{\varepsilon}+\varepsilon, s) \dd s - \int_0^{\hat t_{\varepsilon}}\kappa(\hat t_{\varepsilon}, s) \dd s 
        \Big| \\
        &\hspace{50pt}+ 
        \Big|
        \int_0^{\hat t_{\varepsilon}}\kappa(\hat t_{\varepsilon}, s) \dd s - \int_0^{\hat t_{\varepsilon}}\kappa(\hat t_{\varepsilon}+\varepsilon, s) \dd s
        \Big| \\
        &< \eta + \int_0^\infty \big| \kappa(\hat t_{\varepsilon} + \varepsilon, s) - \kappa(\hat t_{\varepsilon}, s)\big| \dd s < 2 \eta.
    \end{align*}
    Since $\eta$ was arbitrary, this completes the proof.
\end{proof}

Next, we state our main result for deterministic Volterra equations. Using~\cite[Theorem~13.4.7]{gripenberg_londen_staffans_1990} we state a comparison result and using~\cite[Corollary 9.5.6]{gripenberg_londen_staffans_1990} we argue that a given solution must be unique if the vector field is locally Lipschitz continuous. Note that this deterministic result 
can be applied pathwise to the SVEs we are interested in.

\begin{theorem} \label{thm:VIE_comparison}
    Let $\kappa$ be a Volterra kernel of continuous type, let $f \in C([0,\infty), \R^m)$ and $h \in C([0,\infty) \times \R^m, \R^m)$ be locally Lipschitz uniformly in time, i.e.\ for every $n \in \N$  there is some constant $L(n) > 0$ with
    \begin{align*} 
        \sup_{t \in [0,n]} |h(t, x) - h(t,y)| \leq L(n) |x - y|, \quad |x|,|y| \leq n.
    \end{align*} Suppose the Volterra equation
    \begin{align} \label{eq:VE}
        y(t)
        =
        f(t)
        +
        \int_0^t
        \kappa(t,s)h(s,y(s)) \dd s, \quad t \geq 0,
    \end{align}
    admits a solution $y$ up to a time of explosion $T_\infty$. Then the solution $y$ is unique on $[0,T_\infty)$. If $h$ is additionally (component-wise) increasing in the second argument, for a continuous function~$u \in C([0,T_\infty), \R^m)$ with
    \begin{align*} 
		u(t) \overset{(\leq)}{\geq} f(t) + \int_0^t k(t,s) h(s, u(s)) \dd s, \quad t \in [0,T_\infty)
	\end{align*}	
    we have $u(t) \overset{(\leq)}{\geq} y(t)$ on $[0,T_\infty)$. (Here the inequalities are to be understood component-wise.)
\end{theorem}

\begin{proof}
    We argue first that the comparison result holds if the solution $y$ is unique. It suffices to show that~\cite[Theorem~13.4.7]{gripenberg_londen_staffans_1990} is applicable, i.e., the assumptions of \cite[Theorem~12.2.7]{gripenberg_londen_staffans_1990} need to be checked.  By continuity of~$f$ and~$h$, assumptions~(i) and~(ii) are satisfied. For (iii) we note that for any $T > 0$ and $n \in \N$ there exists $C > 0$ such that 
    \begin{align*}
        |\kappa(T,s)h(s,y)|
        \leq
        \kappa(T,s) \sup_{s \in [0,T], |y|\leq n}|h(s,y)|
        \leq C  \kappa(T,s)
    \end{align*}
    for each $s \in [0,T]$ and $|y| \leq n$, and, by definition of a Volterra kernel of continuous type, $\kappa(T,\cdot)$ is integrable on~$[0,T]$. To verify~(iv) we check that the set
    \begin{align*}
        \Big\{
            [0,T] \ni t \mapsto \int_0^t \kappa(t,s)h(s,\psi(s)) \dd s
            \Big| \psi\in C_{f(0),n}([0,T],\R^m)
        \Big\}
    \end{align*}
    is equi-continuous, where $C_{f(0),n}([0,T],\R^m)$ is the set of $\R^m$-valued continuous functions starting in $f(0)$ and modulus bounded by some $n \in \N$. Consider $u, t \in [0,T]$ with $u<t$ and $\psi \in C_{f(0),n}([0,T],\R^m)$. We then have
    \begin{align*}
        \Big|\int_0^t \kappa(t,s)&h(s,\psi(s)) \dd s - \int_0^u \kappa(u,s)h(s,\psi(s)) \dd s\Big| \\
        &\leq 
		\Big|\int_0^u (k(t,s) - k(u,s))h(s, \psi(s)) \dd s + \int_u^t k(t,s) h(s, \psi(s)) \dd s\Big| \\
		&\leq C \left( 
		\int_0^\infty |k(t,s) - k(u,s)| \dd s + \int_u^t |k(t,s)| \dd s\right),
	\end{align*}
    for some $C > 0$. The last two terms converge to zero as $u$ tends to $t$, since $\kappa$ is of continuous type. A similar argument applies if $u > t$, hence we obtain equi-continuity. Thus the comparison result holds.

    We next argue that the solution is unique on $[0,T]$ for $T \in (0, T_\infty)$. Let~$y'$ be another solution defined on $[0,T]$. Note that the function $h$ restricted to the compact set $[0,T] \times [-(\|y\|_\infty + \|y'\|_\infty),\|y\|_\infty + \|y'\|_\infty]^m$ is Lipschitz continuous with constant $L>0$. Thus for $z = |y-y'|$ it holds
    \begin{align*}
        z(t) 
        \leq
        L\int_0^tk(t,s) z(s) \dd s, \quad t \in [0,T].
    \end{align*}
    Knowing that the linear Volterra equation $z(t) = L\int_0^tk(t,s) z(s) \dd s$, $t \in [0,T]$, has zero as its unique solution, see~\cite[Corollary 9.5.6]{gripenberg_londen_staffans_1990}, we can apply the previously established comparison. Thus $z \equiv 0$ and we have shown uniqueness, which completes the proof.
\end{proof}

In our SVEs of interest we allow for forcing terms given by continuous stochastic processes $Z$. We show that for locally Lipschitz vector fields, solutions always exist up to a time of explosion. A path-wise application of the previous result immediately yields path-wise uniqueness for these SVEs. Using~\cite{10.1214/EJP.v12-431} we have uniqueness in law as well.

\begin{theorem}\label{thm:svie}
    Let $\kappa$ be a Volterra kernel of continuous type and consider the filtered probability space $(\Omega, \F, (\F_t)_{t \geq 0}, \P)$. 
    Let $Z$ be a continuous and $\R^m$-valued adapted process starting in zero.
    Further, let $h \in C([0,\infty) \times \R^m,\R^m)$ be locally Lipschitz uniformly in time, see Theorem~\ref{thm:VIE_comparison}.
    Then  the stochastic Volterra integral equation
    \begin{align*}
		Y_t = Z_t + \int_0^t \kappa(t,s) h(s,Y_s) \dd s, \quad  t \geq 0,
	\end{align*}
   admits a path-wise unique solution $Y$, which is adapted and continuous and exists up to a time of explosion $\tau_\infty = \lim_{n \uparrow \infty} \tau_n$, where $\tau_n = \inf\{t \geq 0: |Y_t| > n \}$. Furthermore, this solution is unique in law.
\end{theorem}

\begin{proof}[Proof]
Fix $n \in \N$ for the moment.
    To truncate the function $h$  define 
	\begin{align*}
		h_n(s, y) =
		\begin{cases}
			h(s \land n,y),& \quad \text{for} \quad |y| \leq n, \\
			h\left(s\land n,\frac{y}{|y|}n\right),& \quad \text{for} \quad |y| > n,
		\end{cases}
	\end{align*}
	for $s \geq 0$. The resulting function~$h_n$ is continuous and bounded by the continuity of $h$. Further, we have, for each $t \geq 0$ and all $x, y \in \R^m$,
    \begin{align*}
        |h_n(t, x) - h_n(t, y)| \leq L(n) |x - y|.
    \end{align*}
    Define next the stopping time $\sigma_n = \inf \{ t \geq 0 : |Z_t| \geq n \}$. By Lemma~\ref{app:lem1} we can apply \cite[Theorem~3.1]{zhang2008stochastic}, which yields that there exists a unique adapted processes~$Y^n$ such that for Lebesgue almost all $t \geq 0$ we have
    \begin{align*}
        Y_t^n = Z_t^{\sigma_n} + \int_0^t \kappa(t,s) h_n(s,Y_s^n) \dd s.
    \end{align*}
    By the boundedness of~$h_n$, the right-hand side of the above equation is continuous according to 
    Lemma~\ref{app:lem1}. 
    Hence, we may assume that~$Y^n$ is continuous and satisfies the above equation for all $t \geq 0$ and up to indistinguishability.
    To remove the dependence on~$n$, first define the stopping  times
    \[
    \rho_0 = 0, \quad 
    \rho_n = n  \land \sigma_n
    \land 
    \inf\{t \geq 0 : |Y^n_t| \geq n\}, \quad n \geq 1.
    \]
    We now establish two facts: 
	\begin{enumerate}
		\item[(1)] For each $n \in \N$ and $t \in [0, \rho_{n-1}]$ we have $Y^{n-1}_t = Y^n_t$. To see this, note that for $t \in [0, \rho_{n-1}]$ we have $|Y^{n-1}_t| \leq n-1 < n$ and hence 
		\begin{align*}
			Y^{n-1}_t 
			&= Z_t^{\sigma_{n-1}} + \int_0^t \kappa(t,s) h_{n-1}(s, Y^{n-1}_s) \dd s \\
			&= Z_t^{\sigma_n} + \int_0^t \kappa(t,s) h_n(s, Y^{n-1}_s) \dd s.
		\end{align*}
		Uniqueness and continuity therefore yield $Y^{n-1} = Y^n$ on $[0,\rho_{n-1}]$.
		\item[(2)] We have $\rho_{n-1} \leq \rho_n$ for each  $n \in \N$, which follows immediately from (1): On $[0, \rho_{n-1}]$ the process $Y^n$ agrees with $Y^{n-1}$, which implies $|Y^n_t| \leq n-1$ for each $t \in [0, \rho_{n-1}]$. Hence, for $t \in [0, \rho_{n-1}]$ we have $|Y^n_t| < n$.  
	\end{enumerate}
	Now, we define the process $Y$ by 
	\begin{align*}
		Y_t = \sum_{n = 1}^\infty Y^n_t 1_{[\rho_{n-1}, \rho_n)}(t), \qquad t \geq 0.
	\end{align*}
    Then $Y$ is adapted
 because $Y^n$ is  and $\{\rho_{n-1} \leq t < \rho_n\} \in \mathcal{F}_t$ for each $n \in \N$ and $t \geq 0$. Clearly, $Y$ is continuous on $[0, \lim_{n\uparrow \infty} \rho_n)$ by~(1) and~(2), and jumps to zero afterwards.  
 To see that~$Y$ is a solution of the Volterra equation on $[0, \lim_{n\uparrow \infty} \rho_n)$, we note for $t \in [\rho_{n-1}, \rho_n)$, using $Z^{\sigma_n}_t = Z_t$, that
    \begin{align*}
		Y_t 
		= Y^n_t 
		&= Z_t + \int_0^t \kappa(t,s) h_n(s, Y^n_s) \dd s \\
		&= Z_t + \sum_{l=1}^n \int_{\rho_{l-1} \land t}^{\rho_l \land t} \kappa(t,s) h_n(s, Y^n_s) \dd s \\
		&= Z_t + \sum_{l=1}^n \int_{\rho_{l-1} \land t}^{\rho_l \land t} \kappa(t,s) h(s, Y^l_s) \dd s 
		= Z_t + \int_0^t \kappa(t,s) h(s, Y_s) \dd s .
	\end{align*}
    Hence we have established existence, and path-wise uniqueness immediately follows from Theorem~\ref{thm:VIE_comparison}. Uniqueness in law follows from~\cite[Proposition 2.10]{10.1214/EJP.v12-431}.
    This completes the proof.
\end{proof}

\section{Proof of the reverse direction of Theorem~\ref{thm:rBergomi_multi}}
\label{app:sec:proof_thm:rBergomi_multi}

The proof strategy relies on reducing  the multi-factor model to a a one-dimensional Volterra equation, which is studied in \cite[Theorem~1.1(2)]{Ga19}. We shall  make use of Theorems~\ref{thm:VIE_comparison} and~\ref{thm:svie} from Appendix~\ref{app:SVIE_results}.

\begin{proof}[Proof of Theorem~\ref{thm:rBergomi_multi} (Reverse direction)]
We prove the reverse direction of Theorem~\ref{thm:rBergomi_multi} by contraposition. For convenience we define $\rho = (\sigma^S_{(1:m)})^\top$.
By assumption and upon relabelling, we can assume w.l.o.g.\ that we have $\rho^\top L_1^\top \nu_1 > 0$.
Recalling $W_{1:m} = (W^1, \dots, W^m)^\top$, 
by Lévy's characterization of Brownian motion, the process
\begin{align*}
    B^i 
    = 
    \tfrac{1}{| \nu_i^\top L_i |} \nu_i^\top L_i {W}_{1:m}, 
    \quad i=1,\dots,N,
\end{align*}
is a correlated $\R^N$-valued Brownian motion. Let the matrix $\Gamma \in \R^{N \times m}$ be given such that $B = \Gamma W_{1:m}$. For $i=1,\dots,N$, we define the constants $\beta_i = |\nu_i^\top L_i| > 0$ and the continuous functions
\begin{align*}
    g_i : [0,\infty) \ni t \longmapsto \gamma_i \xi_0(t) \exp\big(- 
        \tfrac{\beta_i^2}{2} t^{2\alpha-1}\big) \in (0,\infty).
\end{align*}
Define the process $\overline Y_t = \int_0^t K_{\alpha,1}(t,s)\dd B_s$ and the function 
\begin{align*}
    h(t,y) = \Big(\sum_{i=1}^N g_i(t) \exp(\beta_i y_i)\Big)^{\frac12}, \quad t \geq 0, \quad y \in \R^N,
\end{align*}
such that the variance process $\psi(t,Y)^2 $ is equivalently given as $h(t, \overline Y_t)^2$.
By Novikov's criterion, the sequence of stopping times defined by $\overline \tau_n = \inf\{t \geq 0: \|\overline Y_t\|_\infty > n\}$ for $n \in \N$ is a localizing sequence for $S$. We define the measures $\dd \P'_n = S_T^{\overline\tau_n}/S_0 \dd \P$ for fixed $T > 0$ and $n \in \N$. 
By Girsanov's theorem, the $\R^N$-valued process 
\begin{align*}
    B^n = 
    B
    -
    \Gamma \rho \int_0^{\cdot \land \tau_n} 
    h(s,\overline Y_s) \dd s
\end{align*}
is a correlated $\P'_n$-Brownian motion. Using this, we see that $\overline Y$ solves the stochastic Volterra equation
\begin{align*} 
    \overline{Y}_t 
    = 
    \int_0^t K_{\alpha,1}(t,s) \dd B^n_s
    + 
    \Gamma \rho \int_0^{t \land \tau_n} K_{\alpha,1}(t,s) h(s,\overline Y_s) \dd s, \quad t \geq 0.
\end{align*}
Analogously, the stochastic Volterra equation
\begin{align} \label{app:B:volt_2}
    \widetilde Y_t = \int_0^t K_{\alpha,1}(t,s) \dd B_s + \Gamma \rho \int_0^t K_{\alpha,1}(t,s) h(s,\widetilde Y_s) \dd s, \quad t \geq 0,
\end{align}
admits a unique solution up to a time of explosion $\widetilde \tau_\infty$ by Theorem~\ref{thm:svie} with announcing sequence $\widetilde \tau_n = \inf\{t \geq 0 :\|\widetilde Y_t \|_\infty > n\}$ for $n \in \N$, satisfying $\lim_{n \uparrow \infty} \widetilde\tau_n = \widetilde\tau_\infty$. By uniqueness of solutions in law, the law of~$\overline Y^{\overline \tau_n}$ under~$\P'_n$ agrees with the law of $\widetilde Y^{\widetilde \tau_n}$ under $\P$. Combining this with the classical argument used by \cite{MR1618849}, we arrive at 
\begin{align} \label{eq:rBergomi_multi:eq2}
    S_0 - \E{S_T} 
    = 
    S_0 \lim_{n \uparrow \infty} \P'_n[\overline \tau_n \leq T]
    = 
    S_0 \P[\widetilde \tau_\infty \leq T],
\end{align}
and it suffices to show that $\P[\widetilde \tau_\infty \leq T] > 0$.

We next  reduce the multi-variate setting to the one-dimensional situation. To this end, consider the one-dimensional Volterra equation 
\begin{align*}
    X_t 
    = 
    \overline Y^1_t + (\Gamma \rho)_1 \int_0^t K_{\alpha,1}(t,s) \sqrt{g_1(s)} \exp(\tfrac12 \beta_1 X_s) \dd s, \quad t \geq 0.
\end{align*}
We claim that the explosion time of the process $X$ is greater than or equal to $\widetilde \tau_\infty$. Note that for $t < \widetilde \tau_\infty$ we have by \eqref{app:B:volt_2} 
\begin{align*}
    \widetilde Y^1_t \geq \overline Y^1_t + (\Gamma \rho)_1 \int_0^t K_{\alpha,1}(t,s) \sqrt{g_1(s)} \exp(\tfrac12 \beta_1 \widetilde Y^1_s) \dd s, \quad t \geq 0.
\end{align*}
By the comparison result Theorem \ref{thm:VIE_comparison} we must have $\widetilde Y^1 \geq X$ 
as long as both processes exist. Since $X \geq \overline{Y}^1$ and 
$\overline{Y}$ is finite, the process $X$ cannot explode before $\widetilde Y^1 $ and the claim holds. Thus it suffices to show that the probability of $X$ exploding before time $T$ is strictly positive under $\P$. This follows from the arguments in \cite[Proof of Theorem~1.1(2)]{Ga19}. Thus $\P[\widetilde \tau_\infty \leq T] > 0$ and $S$ is a strict local martingale by \eqref{eq:rBergomi_multi:eq2}, concluding the proof.
\end{proof}

\section*{Acknowledgements}

We are very grateful to two anonymous referees and an associate editor for their constructive comments.
We thank Hardy Hulley for helpful suggestions.

\end{document}